\documentclass[11pt]{amsart}
\usepackage{amssymb}
\usepackage{amsmath}
\usepackage{amsthm}
\usepackage{latexsym}
\usepackage{amsfonts}
\usepackage{graphicx}
\usepackage{graphics}
\usepackage[T1]{pbsi}

\theoremstyle{definition}
\newtheorem*{Proof}{Proof}

\newcommand{\dis}{\displaystyle}
\newcommand{\ra}{\;\rightarrow\;}

\newcommand{\de}{\delta }

\newcommand{\De} {{\varDelta}}
\newcommand{\e}{\varepsilon }

\newcommand{\La} {{\varLambda}}

\newcommand{\la}{\lambda }

\newcommand{\R}{\mathbb{R}}

\newcommand{\N}{\mathbb{N}}

\newcommand{\ssum}{\sum\limits}

\newcommand{\ld}{\ldots}

\newcommand{\hs}{\hfill$\square$}

\begin{document}

\title[A Hardy type inequality]{A Hardy inequality and applications to reverse H\"{o}lder inequalities for weights on $\R$}
\author{Eleftherios N. Nikolidakis}
\footnotetext{\hspace{-0.5cm}2010 {\em Mathematics Subject Classification.} Primary 26D15; Secondary 42B25.} \footnotetext{\hspace{-0.5cm}{\em Keywords and phrases.} Hardy inequalities, Reverse H\"{o}lder inequalities, weights.}
\date{}
\maketitle
\noindent
{\bf Abstract.} We prove a sharp integral inequality valid for non-negative functions defined on $[0,1]$, with given $L^1$ norm. This is in fact a generalization of the well known integral Hardy inequality. We prove it as a consequence of the respective weighted discrete analogue inequality which proof is presented in this paper. As an application we find the exact best possible range of $p>q$ such that any non-increasing $g$ which satisfies a reverse H\"{o}lder inequality with exponent $q$ and constant $c$ upon the subintervals of $(0,1]$, should additionally satisfy a reverse H\"{o}lder inequality with exponent $p$ and a different in general constant $c'$. The result has been treated in \cite{1} but here we give an alternative proof based on the above mentioned inequality.
\section{Introduction}\label{sec1}
\noindent

During his efforts to simplify the proof of Hilbert's double series theorem, G. H. Hardy \cite{7}, first proved in 1920 the most famous inequality which is known in the literature as Hardy's inequality (see also \cite{10}, Theorem 3.5). This is stated as\vspace*{0.2cm} \\
\noindent
{\bf Theorem A.} {\em If $p>1$, $a_n>0$, and $A_n=a_1+a_2+\cdots+a_n$, $n\in\N$, then
\begin{eqnarray}
\sum^\infty_{n=1}\bigg(\frac{A_n}{n}\bigg)^p<\bigg(\frac{p}{p-1}\bigg)^p\sum^\infty_{n=1}a^p_n. \label{eq1.1}
\end{eqnarray}
Moreover, inequality (\ref{eq1.1}) is best possible, that is the constant and the right side cannot be decreased}.

In 1926, E.Copson, generalized in \cite{3} Theorem A by replacing the arithmetic mean of a sequence by a weighted arithmetic mean. More precisely he proved the following\vspace*{0.2cm} \\
\noindent
{\bf Theorem B.} {\em Let $p>1$, $a_n,\la_n>0$, for $n=1,2,\ld\;.$

Further suppose that $\La_n=\ssum^n_{i=1}\la_i$ and $A_n=\ssum^n_{i=1}\la_ia_i$. Then
\begin{eqnarray}
\sum^\infty_{n=1}\la_n\bigg(\frac{A_n}{\La_n}\bigg)^p\le\bigg(\frac{p}{p-1}\bigg)^p\sum^\infty_{n=1}\la_na^p_n, \label{eq1.2}
\end{eqnarray}
where the constant involved in (\ref{eq1.2}) is best possible.}

In \cite{3}, Copson proves also a second weighted inequality, which as Hardy noted in \cite{8}, can be derived from Theorem B.
From then and until now there have been given several generalizations of the above two inequalities. The first one is given by Hardy and Littlewood who generalized in a specific direction Theorem 1.2 (see \cite{9}). This was generalized further by Leindler in \cite{14}, and by Nemeth in \cite{17}. Also in \cite{16} one can see further generalizations of Hardy's and Copson's series inequalities by replacing means by more general linear transforms. For the study of Copson's inequality one can also see \cite{4}. Additionally, in \cite{5}, Elliot has already proved inequality (\ref{eq1.2}) by similar methods to those that appear in \cite{3}.

There is a continued analogue of Theorem 1.1 (see \cite{10}) which can be stated as \vspace*{0.2cm} \\
\noindent
{\bf Theorem C.} {\em If $p>1$, $f(x)\ge0$ for $x\in [0,+\infty)$ then
\begin{eqnarray}
\int^\infty_0\bigg(\frac{1}{x}\int^x_0f(t)dt\bigg)^pdx<\bigg(\frac{p}{p-1}\bigg)^p\int^\infty_0f^p(x)dx, \label{eq1.3}
\end{eqnarray}
}
Further generalizations of (\ref{eq1.3}) can be seen in \cite{8}.
Other authors have also studied these inequalities in more general forms as it may be seen in \cite{15} and \cite{20}. E. Landau has also studied the above inequality and his work appears in \cite{13}. For a complete discussion of the topic one can consult \cite{12} and \cite{19}.
In this paper we generalize (\ref{eq1.3}) by proving the following \vspace*{0.2cm} \\
\noindent
{\bf Theorem 1.} {\em Let $g:[0,1]\ra\R^+$ be integrable function, $p>1$, and additionally assume that $\int_0^1g=f$. Then the following inequality is true, for any $q$ such that $1\le q\le p$
\begin{eqnarray}
\int^1_0\bigg(\frac{1}{t}\int^t_0g\bigg)^pdt<\bigg(\frac{p}{p-1}\bigg)^q\int^1_0
\bigg(\frac{1}{t}\int^t_0g\bigg)^{p-q}g^q(t)dt-\frac{q}{p-1}f^p  \label{eq1.4}
\end{eqnarray}
Moreover, inequality (\ref{eq1.4}) is sharp in the sense that, the constant $(\frac{p}{p-1})^q$ cannot be decreased, while the constant
$\frac{q}{p-1}$ cannot be increased for any fixed $f$.
}

In fact we are going to prove, an even more general inequality which is the discrete analogue of (\ref{eq1.4}) for the case $q=1$,  which is weighted. This is a generalization of (\ref{eq1.2}) and is described in the following

\noindent
{\bf Theorem 2.} {\em Let $(a_n)_n$ be a sequence of non-negative real numbers. We define for every sequence $(\la_n)_n$ of positive numbers the following quantities $A_n=\la_1a_1+\cdots+\la_na_n$ and $\La_n=\la_1+\cdots+\la_n$. Then the following inequality is true:
\begin{eqnarray}
\sum^N_{n=1}\la_n\bigg(\frac{A_n}{\La_n}\bigg)^p\le\bigg(\frac{p}{p-1}\bigg)
\sum^N_{n=1}\la_na_n\bigg(\frac{A_n}{\La_n}\bigg)^{p-1}-\frac{1}{p-1}\La_N\bigg(\frac{A_N}{\La_N}\bigg)^p,  \label{eq1.5}
\end{eqnarray}
for any $N\in\mathbb{N}$.
}

It is obvious that by setting $\la_n=1$ for every $n\in\N$, in Theorem 2, we reach, for $q=1$, to the discrete analogue of (\ref{eq1.4}), thus generalizing (\ref{eq1.1}) and (\ref{eq1.2}). Then Theorem 1 is an easy consequence, for the case $q=1$, by the use of a standard approximation argument, of $L^1$ functions on $(0,1]$, by simple functions . We then use this result (as can be seen in the sequel) in an effective way to provide a proof of Theorem 1, for any $q\in[1,p]$. We mention also that the opposite problem for negative exponents is treated in \cite{18}.

We believe that Theorem 1 has many applications in many fields and especially in the theory of weights. Our intention in this paper is to describe one of them. We mention the related details.
Let $Q_0\subseteq\R^N$ be a given cube.
Let also $p>1$ and $h:Q_0\rightarrow\R^+$ be such that $h\in L^p(Q_0)$. Then, as it is well known, the following, named as H\"{o}lder's inequality is satisfied
$$
\bigg(\frac{1}{|Q|}\int_Qh\bigg)^p\le\frac{1}{|Q|}\int_Qh^p, \ \ \text{for any cube} \ \ Q\subseteq Q_0.
$$
In this paper we are interested for functions that satisfy a reverse H\"{o}lder inequality. More precisely we say that $h$ satisfies the reverse H\"{o}lder inequality with exponent $q>1$ and constant $c\ge1$ if the following holds
\begin{eqnarray}
\frac{1}{|Q|}\int_Qh^q\le c\cdot\bigg(\frac{1}{|Q|}\int_Qh\bigg)^q \ \ \text{for every cube} \ \ Q\subseteq Q_0.  \label{eq1.6}
\end{eqnarray}
Now in \cite{6} it is proved the following.\vspace*{0.2cm} \\
\noindent
{\bf Theorem A.} {\em Let $1<q<\infty$ and $h:Q_0\rightarrow\R^+$ such that (\ref{eq1.6}) holds. Then there exists $\e=\e(N,q,c)$ such that $h\in L^p$ for any $p$ such that $p\in[q,q+\e)$. Moreover the following inequality holds
$$
\frac{1}{|Q|}\int_Qh^p\le c'\bigg(\frac{1}{|Q|}\int_Qh\bigg)^p,
$$
for any cube $Q\subseteq Q_0$, $p\in[q,q+\e)$ and some constant $c'=c'(N,p,q,c)$}.\medskip

As a consequence the following question naturally arises and is posed in \cite{2} . What is the best possible value of $\e$ ? The problem for the case $N=1$ was solved in \cite{1} for non-increasing functions $g$ and was completed for arbitrary functions in \cite{11}. More precisely in \cite{1} it is shown the following\vspace*{0.2cm} \\
\noindent
{\bf Theorem B.} {\em Let $g:(0,1]\rightarrow\R^+$ be non-increasing which satisfies the following inequality
\begin{eqnarray}
\frac{1}{b-a}\int^b_ag^q\le c\bigg(\frac{1}{b-a}\int^b_ag\bigg)^q,  \label{eq1.7}
\end{eqnarray}
for every $(a,b)\subseteq(0,1]$, where $q>1$ is fixed, and $c$ independent of $a,b$. If we define $p_0>q$ as the root of the following equation
\begin{eqnarray}
\frac{p_0-q}{p_0}\cdot\bigg(\frac{p_0}{p_0-1}\bigg)^q\cdot c=1,  \label{eq1.8}
\end{eqnarray}
we have that $g\in L^p((0,1])$ and $g$ satisfies a reverse H\"{o}lder inequality with exponent $p$, for every $p$ such that $p\in[q,p_0)$. Moreover the result is sharp, that is the value of $p_0$ cannot be increased.}

The problem was solved completely in \cite{11} where the notion of the non-increasing rearrangement of $h$ was used and which is defined as follows:
$$
h^\ast(t)=\sup_{e\subseteq(0,1]\atop|e|\ge t}\Big[\inf_{x\in e}h(x)\Big].
$$
More precisely the following appears in \cite{11}.\vspace*{0.2cm} \\
\noindent
{\bf Theorem C.} {\em Let $h:(0,1]\rightarrow\R^+$, that it satisfies (\ref{eq1.7}), for every $(a,b)\subseteq [0,1]$.
with $q>1$ fixed and $c\ge1$. Then the same inequality is true if we replace $h$ by it's non-increasing rearrangement}.\medskip

It is immediate now that Theorem B and C answer the question as it was posed in \cite{2}, for the case $N=1$.

Our aim in this paper is to give an alternative proof of Theorem B by using Theorem 1. We will prove the following variant of Theorem B which we state as \vspace*{0.2cm} \\
\noindent
{\bf Theorem 3.} {\em Let $g:(0,1]\rightarrow\R^+$ be non-increasing satisfying a reverse H\"{o}lder inequality with exponent $q>1$ and constant $c\ge1$ upon all intervals of the form $(0,t]$. That is the following hold:
\begin{eqnarray}
\frac{1}{t}\int^t_0g^q\le c\cdot\bigg(\frac{1}{t}\int^t_0g\bigg)^q,  \label{eq1.9}
\end{eqnarray}
for any $t\in(0,1]$. Then for every $p\in[q,p_0)$ the following inequality true
\begin{eqnarray}
\frac{1}{t}\int^t_0g^p\le c'\bigg(\frac{1}{t}\int^t_0g\bigg)^p,  \label{eq1.10}
\end{eqnarray}
for any $t\in(0,1]$ where $c'=c'(p,q,c)$ and $p_0$ is defined by (\ref{eq1.5}). As a consequence $g\in L^p$ for every $p\in[q,p_0)$}.\medskip

By the same reasoning we can prove the analogue of Theorem 3, for intervals of the form $(t,1]$. Ending this discussion we mention that in \cite{11} it is proved the following:\vspace*{0.2cm} \\
\noindent
{\bf Theorem D.} {\em Let $g:(0,1]\rightarrow\R^+$ be non-increasing. Then (\ref{eq1.7}) is satisfied for all subintervals of $(0,1]$ iff it is satisfied for all subintervals of the form $(0,t]$ and $[t,1]$}.\medskip

By the above results we conclude that Theorem 3, and its analogue for the intervals of the form $(t,1]$, imply Theorem B.
\section{The Hardy type inequality}\label{sec2}
\noindent

We first present the following which can be seen in \cite{3}.\vspace*{0.2cm} \\
\noindent
{\bf Proof of Theorem 2.}
For each $n\in\N$  define
$$
\De_n=\la_n\bigg(\frac{A_n}{\La_n}\bigg)^p-\frac{p}{p-1}\la_n\bigg(\frac{A_n}{\La_n}\bigg)^{p-1}
a_n=\la_n\De'_n,
$$
where
\[
\De'_n=\bigg(\frac{A_n}{\La_n}\bigg)^p-\frac{p}{p-1}\bigg(\frac{A_n}{\La_n}\bigg)^{p-1}a_n.
\]
Obviously, $a_n=\dfrac{A_n-A_{n-1}}{\la_n}$ for every $n\in\N$, so we have
\begin{align}
\De'_n&=\bigg(\frac{A_n}{\La_n}\bigg)^p-\frac{p}{p-1}\bigg(\frac{A_n}{\La_n}\bigg)^{p-1}
\frac{A_n-A_{n-1}}{\la_n} \nonumber\\
&=\bigg(\frac{A_n}{\La_n}\bigg)^p-\frac{p}{p-1}\bigg(\frac{A_n}{\La_n}\bigg)^p
\frac{\La_n}{\la_n}+\frac{p}{p-1}\bigg(\frac{A_n}{\La_n}\bigg)^{p-1}\frac{A_{n-1}}{\la_n} \nonumber \\
&=\bigg(\frac{A_n}{\La_n}\bigg)^p\bigg[1-\frac{p}{p-1}\cdot\frac{\La_n}{\la_n}\bigg]+
\frac{1}{p-1}\bigg\{p\cdot\bigg(\frac{A_n}{\La_n}\bigg)^{p-1}\frac{A_{n-1}}{\La_{n-1}}
\bigg\}\frac{\La_{n-1}}{\la_n}. \label{eq2.1}
\end{align}
We now use the following elementary inequality
$$
px^{p-1}y\le(p-1)x^p+y^p,
$$
which holds for any $p>1$ and $x,y\ge0$.

We apply it for $x=\dfrac{A_n}{\La_n}$, $y=\dfrac{A_{n-1}}{\La_{n-1}}$, so using (\ref{eq2.1}) we have that:
\begin{align}
\De'_n&\le\bigg(\frac{A_n}{\La_n}\bigg)^p\bigg[1-\frac{p}{p-1}\frac{\La_n}{\la_n}\bigg]+
\frac{1}{p-1}\bigg[(p-1)\bigg(\frac{A_n}{\La_n}\bigg)^p+\bigg(\frac{A_{n-1}}{\La_{n-1}}\bigg)^p
\bigg]\cdot\frac{\La_{n-1}}{\la_n} \nonumber \\
&=\bigg(\frac{A_n}{\La_n}\bigg)^p\bigg[1-\frac{p-1}{p}\frac{\La_n}{\la_n}+\frac{\La_{n-1}}{\la_n}
\bigg]+\frac{1}{p-1}\bigg(\frac{A_{n-1}}{\La_{n-1}}\bigg)^p\frac{\La_{n-1}}{\la_n} \nonumber \\
&=-\frac{1}{p-1}\cdot\frac{\La_n}{\la_n}\bigg(\frac{A_n}{\La_n}\bigg)^p+\frac{1}{p-1}
\frac{\La_{n-1}}{\la_n}\bigg(\frac{A_{n-1}}{\La_{n-1}}\bigg)^p.  \label{eq2.2}
\end{align}
Thus from (\ref{eq2.2}) and the definition of $\De_n$ we conclude
\begin{eqnarray}
\De_n\le\frac{1}{p-1}\La_{n-1}\bigg(\frac{A_{n-1}}{\La_{n-1}}\bigg)^p-\frac{1}{p-1}
\La_n\bigg(\frac{A_n}{\La_n}\bigg)^p,  \label{eq2.3}
\end{eqnarray}
This holds for every $n\in\N$, $n\ge2$.

It is immediate now that for $n=1$ we have the following equality
\begin{eqnarray}
\De_1=-\frac{1}{p-1}\La_1\bigg(\frac{A_1}{\La_1}\bigg)^p.  \label{eq2.4}
\end{eqnarray}
For any $N\in\N$ we sum (\ref{eq2.3}) from $n=2$ to $N$ and add also the equality (\ref{eq2.4}), so we conclude after making the appropriate
cancellations, inequality (\ref{eq1.5}) of Theorem 2.  \hs

The following is an easy consequence of the above result

{\bf Corollary 1:} {\em
Let $g:[0,1]\ra\R^+$ be integrable function, $p>1$ and additionally assume that $\int_0^1g=f$. Then the following inequality is true
\begin{eqnarray}
\int^1_0\bigg(\frac{1}{t}\int^t_0g\bigg)^pdt\le\bigg(\frac{p}{p-1}\bigg)\int^1_0
\bigg(\frac{1}{t}\int^t_0g\bigg)^{p-1}g(t)dt-\frac{1}{p-1}f^p.  \label{eq2.5}
\end{eqnarray}
}
We proceed now to the

{\bf Proof of Theorem 1.}
For any $s\in[0,p]$ we define by $I_s$ by
$$I_s=\int^1_0\bigg(\frac{1}{t}\int^t_0g\bigg)^{p-s}g^s(t)dt,$$
for any $g:[0,1]\ra\R^+$ integrable function, such that $\int_0^1g=f$. Then, for the proof of inequality (\ref{eq1.4}), we just need to prove that
$$I_0\le\bigg(\frac{p}{p-1}\bigg)^qI_q-\frac{q}{p-1}f^p,$$
for any $q\in(1,p]$.

We write $$I_1=\int^1_0g(t)\bigg(\frac{1}{t}\int^t_0g\bigg)^{(p-q)/q}\bigg(\frac{1}{t}\int^t_0g\bigg)^{p-\frac{p}{q}}dt.$$

We then apply in the above integral H\"{o}lder's inequality, with exponents $q,\frac{q}{q-1}$, and we have as a consequence that
\begin{eqnarray}
I_1\le I_q^{1/q}I_0^{(q-1)/q}.  \label{eq2.6}
\end{eqnarray}
Additionally from Corollary 1 we obtain
\begin{eqnarray}
I_0\le\frac{p}{p-1}I_1-\frac{1}{p-1}f^p . \label{eq2.7}
\end{eqnarray}
We consider now the difference $L_q=I_0-(\frac{p}{p-1})^qI_q$.
We need to prove that $L_q\le-\frac{q}{p-1}f^p$.
By using the inequalities (\ref{eq2.6}) and (\ref{eq2.7}) we have that
\begin{align}
L_q&\le I_0-\bigg(\frac{p}{p-1}\bigg)^q\frac{I_1^q}{I_0^{q-1}} \nonumber \\
&\le I_0-\bigg(\frac{p}{p-1}\bigg)^q I_0^{-q+1}\bigg(\frac{p-1}{p}I_0+\frac{1}{p}f^p\bigg)^q.        \label{eq2.8}
\end{align}
We define now the following function of the variable $x>0$:
$$G(x)=x-\bigg(\frac{p}{p-1}\bigg)^q x^{-q+1}\bigg(\frac{p-1}{p}x+\frac{1}{p}f^p\bigg)^q.$$
Then
$$G(x)=x-x^{-q+1}\bigg(x+\frac{1}{p-1}f^p\bigg)^q,$$
so that
$$G'(x)=1+(q-1)\bigg(1+\frac{f^p}{(p-1)x}\bigg)^q-q\bigg(1+\frac{f^p}{(p-1)x}\bigg)^{q-1}.$$
Now we consider the following function of the variable $t\ge1$: $F(t)=1+(q-1)t^q-qt^{q-1}$. Then $F'(t)=q(q-1)t^{q-2}(t-1)>0$, for every $t>1$.
Thus $F$ is strictly increasing on its domain, so that $F(t)>F(1)=0$, for any $t>1$. We immediately  conclude that $G'(x)>0$, for every $x>0$.
As a consequence $G$ is strictly increasing on $(0,+\infty)$. We evaluate now $\underset{x\rightarrow +\infty}{lim}G(x)=l$.
We have that
\begin{align}
l&=\underset{x\rightarrow +\infty}{lim} x\bigg[1-\bigg(1+\frac{f^p}{(p-1)x}\bigg)^q\bigg] \nonumber \\
&=\underset{x\rightarrow +\infty}{lim}\frac{1-\bigg(1+\frac{yf^p}{p-1}\bigg)^q}{y}=-\frac{q}{p-1}f^p,         \label{eq2.9}
\end{align}
by using De'l Hospital's rule. Thus since $G$ is strictly increasing on $(0,+\infty)$, we have that $G(x)<-\frac{q}{p-1}f^p$, for any $x>0$.
Thus (\ref{eq2.8}) yields $L_q<-\frac{q}{p-1}f^p$, which is inequality (\ref{eq1.4}). We now prove its sharpness.

We let
\[
J'_0=\int^1_0\bigg(\frac{1}{t}\int^t_0g\bigg)^pdt, \ \ \text{and} \ \ J'_q=
\int^1_0\bigg(\frac{1}{t}\int^t_0g\bigg)^{p-q}g^q(t)dt \ \ \text{for any} \ \ 1\le q\le p.
\]
Let also $g=g_a$, where $g_a$ is defined for any $a\in(0,1/p)$, by $g_a(t)=t^{-a}$, $t\in(0,1]$.

Then for every $t\in(0,1]$ we have that
\[
\frac{1}{t}\int^t_0g_a=\frac{1}{1-a}g_a(t) \ \ \text{and so}
\]
\begin{eqnarray}
\frac{J'_0}{J'_q}=\frac{\Big(\dfrac{1}{1-a}\Big)^p\dis\int^1_0g_a^pdt}
{\Big(\dfrac{1}{1-a}\Big)^{p-q}\dis\int^1_0g^p_adt}=\bigg(\frac{1}{1-a}\bigg)^q. \label{eq2.10}
\end{eqnarray}
Letting $a\rightarrow1/p^{-}$ in (\ref{eq2.10}) we obtain that the constant $(\frac{p}{p-1})^q$, on the right of inequality (\ref{eq1.4}), cannot be
decreased. We now prove the second part of the sharpness of Theorem 1. For this purpose we define for any fixed $f>0$, and any $a\in(0,1/p)$,
the function $g_a(t)=f(1-a)t^{-a}$, for every $t\in(0,1]$. Then it is easy to see that $\int^1_0g_a=f$, $\frac{1}{t}\int^t_0g_a(u)du=\frac{1}{1-a}g_a(t)$, and that $\int^1_0g_a^p=\frac{f^p(1-a)^p}{1-ap}$.
We consider now the difference
\begin{align}
L_q(a)&=\int^1_0\bigg(\frac{1}{t}\int^t_0g_a\bigg)^{p}dt-\bigg(\frac{p}{p-1}\bigg)^q\int^1_0\bigg(\frac{1}{t}\int^t_0g_a\bigg)^{p-q}g_a^q(t)dt
 \nonumber \\
&=\bigg(\frac{1}{1-a}\bigg)^p\int^1_0g_a^p-\bigg(\frac{p}{p-1}\bigg)^q\bigg(\frac{1}{1-a}\bigg)^{p-q}\int^1_0g_a^p   \nonumber \\                                                       &=\frac{\bigg(\frac{1}{1-a}\bigg)^{p-q}f^p\bigg[\bigg(\frac{1}{1-a}\bigg)^{q}-\bigg(\frac{p}{p-1}\bigg)^{q}\bigg]}{1-ap}.   \label{eq2.11}
\end{align}
Letting now $a\rightarrow1/p{-}$, we immediately see, by an application of De'l Hospital's rule that $L_q(a)\rightarrow-\frac{q}{p-1}f^p$.
We have just proved that the constant $\frac{q}{p-1}$, appearing in front of $f^p$, cannot be increased. That is, both constants appearing
on the right of (\ref{eq1.4}) are best possible.     \hs

\section{Applications to reverse H\"{o}lder inequalities}\label{sec3}

We will need first a preliminary lemma which in fact holds under some additional hypothesis for $g$ even if it is not decreasing, which can be proved using integration by parts. We present a version that we will need below which is proved by measure integration techniques. More precisely we will prove the following \vspace*{0.2cm} \\
\noindent
{\bf Lemma 1.} {\em Let $g:(0,1]\ra\R^+$ be a non-increasing function. Then the following inequality is true for any $p>1$ and every $\de\in(0,1)$}
\begin{eqnarray}
\int^\de_0\bigg(\frac{1}{t}\int^t_0g\bigg)^pdt=-\frac{1}{p-1}\bigg(\int_0^\de g\bigg)^p\frac{1}{\de^{p-1}}+\frac{p}{p-1}\int^\de_0\bigg(\frac{1}{t}\int^t_0g\bigg)^{p-1}g(t)dt.\hspace*{-1cm} \label{3.1}
\end{eqnarray}
\begin{Proof}
By using Fubini's theorem it is easy to see that
\begin{eqnarray}
\int^\de_0\bigg(\frac{1}{t}\int^t_0g\bigg)^pdt=\int^{+\infty}_{\la=0}p\la^{p-1}\bigg|\bigg\{t\in(0,\de]:
\frac{1}{t}\int^t_0g\ge\la\bigg\}\bigg|dt. \label{3.2}
\end{eqnarray}
Let now $\dfrac{1}{\de}\dis\int^\de_0g=f_\de\ge f=\dis\int^1_0g$. Then
\[
\begin{array}{l}
  \dfrac{1}{t}\dis\int^t_0g>f_\de, \ \ \forall\;t\in(0,\de) \ \ \text{while} \\ [2ex]
  \dfrac{1}{t}\dis\int^t_0g\le f_\de, \ \ \forall\;t\in[\de,1].
\end{array}
\]
Let $\la$ be such that: $0<\la<f_\de$. Then for every $t\in(0,\de]$ we take $\dfrac{1}{t}\dis\int^t_0g\ge\dfrac{1}{\de}\dis\int^\de_0g=f_\de>\la$. Thus
\[
\bigg|\bigg\{t\in(0,\de]:\frac{1}{t}\int^t_0g\ge\la\bigg\}\bigg|=|(0,\de]|=\de.
\]
Now for every $\la>f_\de$ there exists unique $a(\la)\in(0,\de)$ such that $\dfrac{1}{a(\la)}\dis\int^{a(\la}_0g=\la$. It's existence is quaranteeded by the fact that $\la>f_\de$, that $g$ is non-increasing and that $g(0^+)=+\infty$ which may without loss of generality be assumed (otherwise we work for the $\la$'s on the interval $(0,\|g\|_\infty])$.
Then
\[
\bigg\{t\in(0,\de]:\frac{1}{t}\int^t_0g\ge\la\bigg\}=(0,a(\la)].
\]
Thus, from the above and (3.23) we conclude that
\begin{align}
\int^\de_0\bigg(\frac{1}{t}\int^t_0g\bigg)^pdt&=\int^{f_\de}_{\la=0}p\la^{p-1}\cdot
\de\cdot d\la+\int^{+\infty}_{\la=f_\de}p\la^{p-1}a(\la)d\la \nonumber\\
&=\de(f_\de)^p+\int^{+\infty}_{\la=f_\de}p\la^{p-1}\frac{1}{\la}\bigg(\int^{a(\la)}_0
g(u)du\bigg)d\la                           \label{3.3}
\end{align}
by the definition of $a(\la)$. As a consequence (\ref{eq3.3}) gives
\begin{align*}
\int^\de_0\bigg(\frac{1}{t}\int^t_0g\bigg)^pdt&=\frac{1}{\de^{p-1}}\bigg(\int^\de_0g\bigg)^p
+\int^{+\infty}_{\la=f_\de}p\la^{p-2}\bigg(\int^{a(\la)}_0g(u)du\bigg)d\la \\
&=\frac{1}{\de^{p-1}}\bigg(\int^\de_0g\bigg)^p+\int^{+\infty}_{\la=f_\de}p\la^{p-2}
\bigg(\int_{\{u\in(0,\de]:\atop\frac{1}{u}\int^u_0g\ge\la\}}g\bigg)d\la\\
&=\frac{1}{\de^{p-1}}\bigg(\int^\de_0g\bigg)^p+\frac{p}{p-1}\int^\de_0g(t)\Big[
\la^{p-1}\Big]^{\frac{1}{t}\int^t_0g}_{\la=f_\de}dt\\
&=\frac{1}{\de^{p-1}}\bigg(\int^\de_0g\bigg)^p+\frac{p}{p-1}\bigg[\int^\de_0
\bigg(\frac{1}{t}\int^t_0g\bigg)^{p-1}g(t)-\bigg(\int^\de_0g(t)dt\bigg)f_\de^{p-1}\bigg]\\
&=-\frac{1}{p-1}\frac{1}{\de^{p-1}}\bigg(\int^\de_0g\bigg)^p+\frac{p}{p-1}\int^{\de}_0\bigg(
\frac{1}{t}\int^t_0g\bigg)^{p-1}g(t)dt,
\end{align*}
where in the third equality we have used Fubini's theorem and the fact that $\dfrac{1}{\de}\dis\int^\de_0g=f_\de$.
In this way we derived (3.22). \hs
\end{Proof}
We are now able to give the\vspace*{0.2cm} \\
\noindent
{\bf Proof of Theorem 3.} Suppose we are given $g:(0,1]\rightarrow\R^+$ non-increasing and $\de\in(0,1]$. Our hypothesis for $g$ is (\ref{eq1.10}) or that:
$$
\frac{1}{t}\int^t_0g^q\le c\cdot\bigg(\frac{1}{t}\int^t_0g\bigg)^q,  \ \ \text{for every} \ \ t\in(0,1].
$$
Let now $p>q$ and set $a=p/q>1$.

We apply Lemma 1 with $g^q$ in place of $g$ and $a$ in that of $p$. We conclude that:
\[
\int^\de_0\bigg(\frac{1}{t}\int^t_0g^q\bigg)^{p/q}dt\!\le\!-\frac{q}{p-q}\frac{1}{\de^{p/q-1}}
\bigg(\int^\de_0g^q\bigg)^{p/q}
\!+\!\frac{p}{p-q}\int^\de_0\bigg(\frac{1}{t}\int^t_0g^q\bigg)^{p/q-1}
g^q(t)dt \Rightarrow
\]
\begin{eqnarray}
\frac{1}{\de}\int^\de_0\bigg[\bigg(\frac{1}{t}\int^t_0g\bigg)^{p/q-1}g^q(t)-
\frac{p-q}{p}\bigg(\frac{1}{t}\int^t_0g^q\bigg)^{p/q}\bigg]dt\le\frac{q}{p}
\bigg(\frac{1}{\de}\int^\de_0g^q\bigg)^{p/q}. \hspace*{-2cm} \label{eq3.4}
\end{eqnarray}
Define now for every $y>0$ the function $\phi_y$ with variable $x$ by $\phi_y(x)=x^{p/q-1}y-\dfrac{p-q}{p}x^{p/q}$, for $x\ge y$.

Then
\[
\phi'_y(x)=(p/q-1)x^{p/q-2}y-(p/q-1)x^{p/q-1}=(p/q-1)x^{p/q-2}(y-x)\le0, \ \ \text{for} \ \ x\ge y.
\]
Thus
\begin{eqnarray}
y\le x\le z\;\Rightarrow\;\phi_y(x)\ge \phi_y(z)  \label{eq3.5}
\end{eqnarray}
Let us now set in (\ref{eq3.5})
\[
x=\frac{1}{t}\int^t_0g^q, \ \ y=g^q(t), \ \ z=c\bigg(\frac{1}{t}\int^t_0g\bigg)^q \ \ \text{for any} \ \ t\in(0,1].
\]
Then $y\le x\le z{\Rightarrow}$
\begin{align*}
\bigg(\frac{1}{t}\int^t_0g^q\bigg)^{p/q-1}g^q(t)-\frac{p-q}{p}\bigg(\frac{1}{t}\int^t_0
g^q\bigg)^{p/q}\ge&\; c^{p/q-1}\bigg(\frac{1}{t}\int^t_0g\bigg)^{p-q}
g^q(t)\\
&-\frac{p-q}{p}c^{p/q}\bigg(\frac{1}{t}\int^t_0g\bigg)^p, \ \ \forall\;t\in(0,1].
\end{align*}
As a consequence (\ref{eq3.4}) gives, by using the last inequality the following
$$
\frac{1}{\de}\int^\de_0\bigg(\frac{1}{t}\int^t_0g\bigg)^{p-q}g^q(t)dt\le c\cdot
\frac{p-q}{p}\cdot\frac{1}{\de}\int^\de_0\bigg(\frac{1}{t}\int^t_0g\bigg)^pdt+
\frac{q}{p}c\bigg(\frac{1}{\de}\int^\de_0g\bigg)^p, \hspace*{-1.5cm}
$$
We use now the inequality,
\[
\frac{1}{\de}\int^\de_0\bigg(\frac{1}{t}\int^t_0g\bigg)^pdt\le\bigg(\frac{p}{p-1}\bigg)^q
\frac{1}{\de}\int^\de_0\bigg(\frac{1}{t}\int^t_0g\bigg)^{p-q}g^q(t)dt
\]
which is a consequence of Theorem 1.

We conclude that if $p_0$ is defined by (\ref{eq1.8}), for any $p\in[q,p_0)$, the following holds
\[
\bigg[1-c\frac{p-q}{p}\bigg(\frac{p}{p-1}\bigg)^q\bigg]\frac{1}{\de}\int^\de_0
\bigg(\frac{1}{t}\int^t_0g\bigg)^{p-q}g^a(t)dt\le\frac{q}{p}c\bigg(\frac{1}{\de}
\int^\de_0g\bigg)^p,
\]
where $1-c\dfrac{p-q}{p}\Big(\dfrac{p}{p-1}\Big)^q=k_p>0$, for every such $p$. This becomes
\begin{eqnarray}
\frac{1}{\de}\int^\de_0\bigg(\frac{1}{t}\int^t_0g\bigg)^{p-q}g^q(t)dt\le\frac{q\cdot c}{p\cdot k_p}\bigg(\frac{1}{\de}\int^\de_0g\bigg)^p,  \label{eq3.6}\vspace*{-0.2cm}
\end{eqnarray}
for any $\de\in(0,1]$, and any $p\in[q,p_0)$.

On the other hand $\dfrac{1}{t}\dis\int^t_0g\ge g(t)$, since $g$ is non-increasing, thus (\ref{eq3.6}) $\Rightarrow$
\[
\frac{1}{\de}\int^\de_0g^p\le\frac{q\cdot c}{p\cdot k_p}\bigg(\frac{1}{\de}\int^\de_0g\bigg)^p,
\]
for any $\de\in(0,1]$ and $p$ such that $q\le p<p_0$, which is an inequality of the form of (\ref{eq1.10}), for suitable $c'>1$.

So the first part of  Theorem 3 is now proved. We continue with the sharpness of the result. For this reason we define for any fixed $c\ge1$ and $q>1$ the following function $g_a(t)=t^{-a}$ where $a=1/p_0$, where $p_0$ is defined by (\ref{eq1.5}). Then it is easy to see that $\dfrac{1}{t}\dis\int^t_0g^q=c\Big(\dfrac{1}{t}\dis\int^t_0g\Big)^q$, for every $t\in(0,1]$. It is obvious now that $g\notin L^{p_0}((0,1])$.
Thus, $p_0$ cannot be increased and Theorem 2 is proved.  \hs
\vspace*{1cm}
\noindent
Nikolidakis Eleftherios\vspace*{0.1cm}\\
Post-doctoral researcher\vspace*{0.1cm}\\
National and Kapodistrian University of Athens\vspace*{0.1cm}\\
Department of Mathematics\vspace*{0.1cm}\\
Panepistimioupolis, GR 157 84\vspace*{0.1cm}\\
Athens, Greece \vspace*{0.1cm}\\
E-mail address:lefteris@math.uoc.gr


\begin{thebibliography}{99}
%
\bibitem{1} L. D' Appuzzo and C. Sbordone, {\em Reverse H\"{o}lder inequalities. A sharp result \em}. Rendiconti Math., {\bf 10}, Ser. VII, (1990), 357-366.
%
\bibitem{2} B. Bojarski, {\em Remarks on the stability of reverse H\"{o}lder inequalities and quasiconformal mappings}, Ann. Acad. Sci. Fenn. A Math. {\bf10},
(1985), 291-296.
%
\bibitem{3} E. T. Copson, {\em Note on series of positive terms}, J. London Math. Soc. {\bf2} (1927), 9-12 and {\bf 3} (1928), 49-51.
%
\bibitem{4} E. T. Copson, {\em Some integral inequalities}, Proc. Roy. Soc. Edinburgh Sect. A {\bf 75} (1975/1976), 157-164.
%
\bibitem{5} E. B. Elliot, {\em A simple exposition of some recently proved facts as to convergency \em}. J. London Math. Soc., {\bf 1}, (1926), 93-96.
%
\bibitem{6} F. W. Gehring, {\em The $L^p$ inequalities of the partial derivatives of a quasiconformal mapping}, Acta Math., {\bf 130}, (1973), 265-27.
%
\bibitem{7} G. H. Hardy, {\em Note on a theorem of Hilbert}, Math Z. {\bf 6} (1920), 314-317.
%
\bibitem{8} G. H. Hardy, {\em Notes on some points in the integral calculus, L. X. An inequality between integrals \em}. Messenger of Math. {\bf 57}, (1928), 12-16.
%
\bibitem{9} G. H. Hardy and J. E. Littlewood, {\em Elementary theorems concerning power series with positive coefficients and moment constants of positive functions}, J. Reine Angew. Math. {\bf 157} (1927), 141-158.
%
\bibitem{10} G. H. Hardy, J. E. Littlewood and G. Polya, {\em Inequalities \em}, Cambridge University Press, Cambridge (1934).
%
\bibitem{11} A. A. Korenovskii, {\em The exact continuation of a Reverse H\"{o}lder inequality and Muckenhoupt's conditions \em}, Math. Notes {\bf 52}, (1992) 1192-1201.
%
\bibitem{12} A. Kufner, L. Maligranda and L. E. Persson, {\em The prehistory of the Hardy inequality \em}, Amer. Math. Monthly, Vol. {\bf 113}, No 8, (2006), 715-732.
%
\bibitem{13} E. Landau, {\em A note on a theorem concerning series of positive terms: Extract from a letter of Prof. E. Landau to Prof. I. Schur \em}, J. London Math. Soc. {\bf 1}, (1926), 38-39.
%
\bibitem{14} L. Leindler, {\em Generalization of inequalities of Hardy and Littlewood}, Acta Sci. Math {\bf 31} (1970), 279-285.
%
\bibitem{15} N. Levinson, {\em Generalizations of an inequality of Hardy}, Duke Math. J. {\bf 31} (1964), 389-394.
%
\bibitem{16} E. R. Love, {\em Generalizations of Hardy's and Copson's inequalities}, J. London Math. Soc. {\bf30} (1984), 431-440.
%
\bibitem{17} J. Nemeth, {\em Generalizations of the Hardy-Littlewood inequality}, Acta Sci. Math. (Szeged) {\bf 32} (1971), 295-299.
%
\bibitem{18} E. N. Nikolidakis, {\em A sharp integral Hardy type inequality and applications to Muckenhoupt weights on R}, Ann. Acad. Sci. Fenn. Math. Vol 39, (2014), 887-896.
%
\bibitem{19} B. G. Pachpatte, {\em Mathematical Inequalities \em}, North Holland Mathematical Library, Vol {\bf 67}, (2005).
%
\bibitem{20} J. B. G. Pachpatte, {\em On a new class of Hardy type inequalities}, Proc. Roy. Soc. Edinburgh Sect. A {\bf105} (1987), 265-274.
%
\end{thebibliography}
\end{document}